\documentclass[12pt]{amsart}

\numberwithin{equation}{section}
\usepackage{graphicx}
\usepackage{amsfonts}
\usepackage{amssymb}
\usepackage{amsmath}
\usepackage{lmodern}
\usepackage{wrapfig}
\usepackage{placeins}
\usepackage{epigraph}
\usepackage{caption}
\usepackage{hyperref}

\begin{document}
 
 \title{The  Black Rabbits of Fibonacci}
 
\author{ A. Solyanik}

\date{24 December 2024}
\dedicatory{A Christmas gift to my dear friend, Prof. Alexander Stokolos}

\address{Department of Applied Mathematics, \ Odessa National Polytechnic University, \ Odessa, \ Ukraine}
\email{solyanik@opu.ua}

\newtheorem{remark}{Remark}
\newtheorem{theorem}{Theorem}
\newtheorem{lemma}{Lemma}
\newtheorem{lemma*}{Lemma*}
\newtheorem{proposition}{Proposition}
\newtheorem{corollary}{Corollary}
\newtheorem{example}{Example}

\maketitle
\begin{abstract}
In this note, we use a toy problem of detecting cycles of length two in a tent map to highlight some curious phenomena in the behavior of discrete dynamical systems.

This work presents no new results or proofs — only computer experiments and illustrations. Thus, it serves as light reading and does not aim to be a scientific paper but is rather educational in nature. For this reason, it is accompanied by numerous illustrations.

\end{abstract}

\epigraph {Alice was beginning to get very tired of sitting by her sister on the bank, and of having nothing to do: once or twice she had peeped into the book her sister was reading, but it had no pictures or conversations in it, “and what is the use of a book,” thought Alice, “without pictures or conversations?"}{Lewis Carroll, Alices Adventures in Wonderland}

\section{Introduction}

This story began when the author became interested in the efficiency of stabilization algorithms he had recently proposed for stabilizing cycles in discrete dynamical systems.

One of the simplest discrete dynamical systems is generated by the tent map:

\[
T_h(x) = 
\begin{cases} 
hx & \text{if } x \in [0, 0.5], \\
h(1 - x) & \text{if } x \in (0.5, 1],
\end{cases}
\]
with a parameter $h\in (1, 2]$.

Professor Alexander Stokolos drew the author's attention to this dynamical system, which was subsequently used as a benchmark.

Unexpectedly, the author discovered that this system exhibits rather intriguing behavior. After spending a couple of days exploring it, the author uncovered several facts that readers can easily verify and prove.

As an introduction to the theory of discrete dynamical systems, we highly recommend the excellent books \cite{cg}, \cite{ce} and \cite{mt}

\begin{figure}[htbp]
    \centering
    \includegraphics[width=\textwidth]{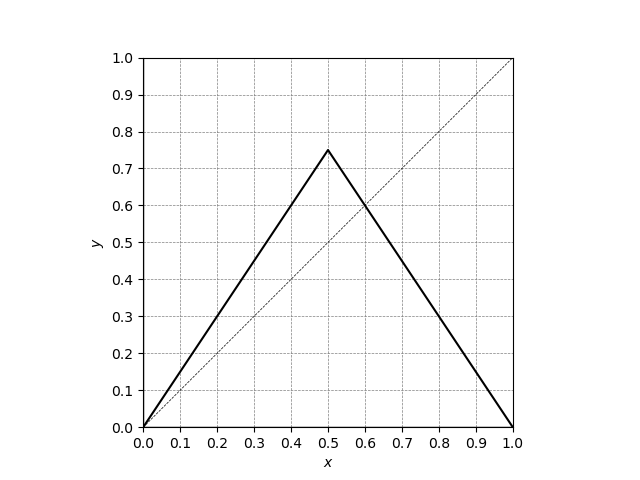} 
    \caption{Plot of $T_h(x)$ for $h=1.5$ and line $y=x$}
    \label{fig:plot}
\end{figure}

\FloatBarrier 
The tent map, for every \( h > 1 \), has one fixed point and also possesses cycles of lengths \( 2, 2^2, 2^3, 2^4, \ldots \). All these cycles appear simultaneously as \( h \) passes the point \( 1 \). It is also quite unexpected that the general formulas for the points of 2- and 4-cycles are simple rational functions, whereas the general formula for an 8-cycle is represented by a rational spline.

We recall Sharkovsky's theorem: let the natural numbers be ordered as follows:
$$
3 \prec 5 \prec 7 \prec \cdots \prec 2 \cdot 3 \prec 2 \cdot 5 \prec 2 \cdot 7 \prec \cdots \prec 2^2 \cdot 3 \prec 2^2 \cdot 5 \prec 2^2 \cdot 7 \prec \cdots \prec 2^3 \prec 2^2 \prec 2 \prec 1
$$

If a continuous map \( f \) from \([0,1]\) to itself has a cycle of period \( p \), and \( p \prec q \), then \( f \) also has a cycle of period \( q \).

The first cycle of length \( 3 \) appears when \( h \) equals the golden ratio, \( \varphi = \frac{\sqrt{5}+1}{2} \approx 1.61803398\ldots \). The first cycle of period \( 5 \) for the tent map emerges when \( h \) equals the largest positive root of the equation \( h^4 - h^3 - h^2 + h - 1 = 0 \), which is approximately \( 1.51287639\ldots \). The first cycle of period \( 7 \) appears when \( h \) equals the largest positive root of the equation \( h^6 - h^5 - h^4 + h^3 - h^2 + h - 1 = 0 \), approximately \( 1.46557123\ldots \), and so on. A cycle of period \( 6 \) appears when \( h = \sqrt{\varphi} \approx 1.27019639\ldots \), and so on.

All these cycles are unstable, which makes this dynamical system an excellent benchmark for testing stabilization algorithms.

The motivation for this investigation stemmed from the challenge of detecting cycles in chaotic time series.

Consider an observed, unpredictable, and chaotic time series that appears as follows:

\begin{figure}[htbp]
    \centering
    \includegraphics[width=\textwidth]{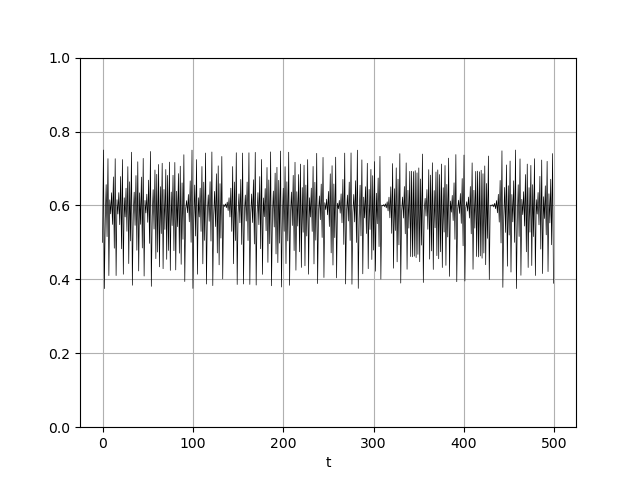} 
    \caption{Observed time series}
    \label{fig_chaotic}
\end{figure}

\FloatBarrier

It is assumed that this time series results from measurements of some observable within a real-world system — be it physical, ecological, economic, financial, biological, etc. Since we are unaware of the governing laws, we consider the system as a black box.

Is it possible to stabilize the system's state using a stabilization algorithm? Can we determine whether the system has (unstable) cycles and their order? Can we find all of these cycles (at least theoretically)? Is it feasible to ascertain this solely by analyzing the time series plots without influencing the system?

It turns out that affirmative answers can be provided to all these questions.

Worth noting, the time series  in Figure \ref{fig_chaotic} is simply the sequence \( T^t_h(0.5) \) for \( h = 1.5 \) and natural numbers \( t = 0, 1, 2, \ldots \).

\section{Detecting Cycles}

The starting point for evaluating the effectiveness of stabilization algorithms was the stabilization of a cycle of length 2 for the map $T_{1.5}$.

For every \(1 < h \leq 2\), there are three fixed points of \(T_h^2\): one is \(\frac{h}{h+1}\), a fixed point of \(T_h\), and two points from the cycle, which are \(\frac{h}{1+h^2}\) and \(\frac{h^2}{1+h^2}\). For the parameter value \(h=1.5\), they are: fixed point \(\frac{3}{5} = 0.6\) and cycle, consisting of two points \(\frac{6}{13}= 0.\overline{461538}\) and \(\frac{9}{13}=0.\overline{692307}\). 

Hereafter, to avoid carrying the subscript, we will denote \( T_h \) with \( h = 1.5 \) simply as \( T \).

\begin{figure}[htbp]
    \centering
    \includegraphics[width=\textwidth]{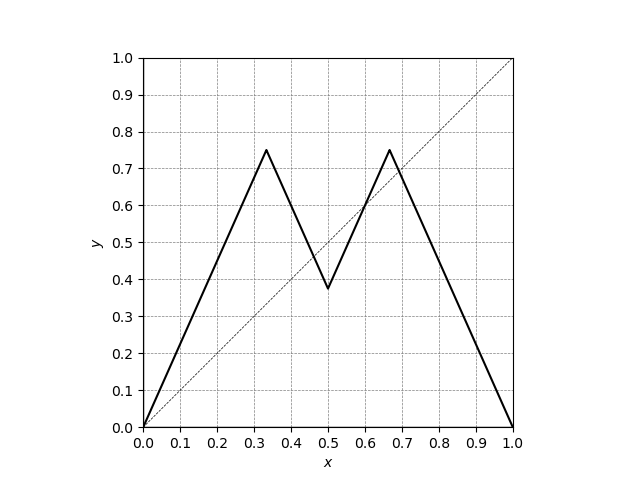} 
    \caption{Plot of $T^2(x)$ and line $y=x$}
    \label{fig:plot}
\end{figure}

\FloatBarrier 

The cycle \(6/13 \rightarrow 9/13 \rightarrow 6/13 \) is unstable. While it is straightforward to verify that \( T \left( 6/13 \right) =9/13 \) and \( T\left(9/13\right) = 6/13 \), any point \( x \) arbitrarily close to \( 6/13 \) will have its trajectory \( T^n(x) \) initially exhibit cyclic behavior but eventually transition to chaotic dynamics.

The following two plots illustrate the trajectories \( T^n(x) \) and \( (T^2)^n(x) \) for \( n \leq 50 \). The starting point \( x = 0.461538 \) was chosen within an interval of radius \( 10^{-4} \) around \( 6/13 \), which is a fixed point of \( T^2 \) and one point of the 2-cycle for \( T \).

\begin{figure}[htbp]
    \centering
    \begin{minipage}[t]{0.48\textwidth}
        \centering
        \includegraphics[width=\textwidth]{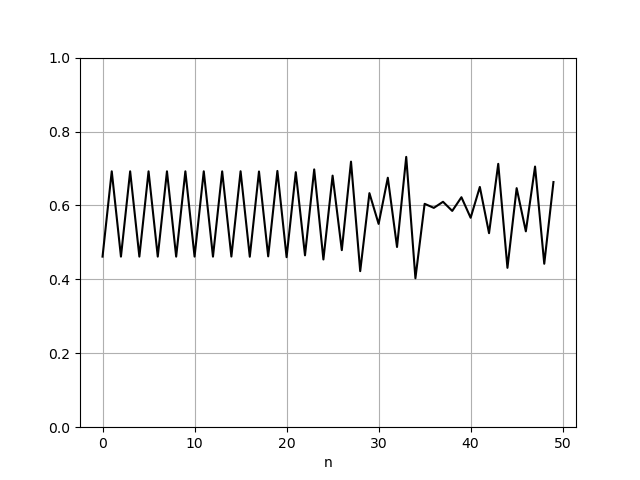}
        \caption*{ $T^n(x_0)$}
        \label{fig:image1}
    \end{minipage}
    \hfill
    \begin{minipage}[t]{0.48\textwidth}
        \centering
        \includegraphics[width=\textwidth]{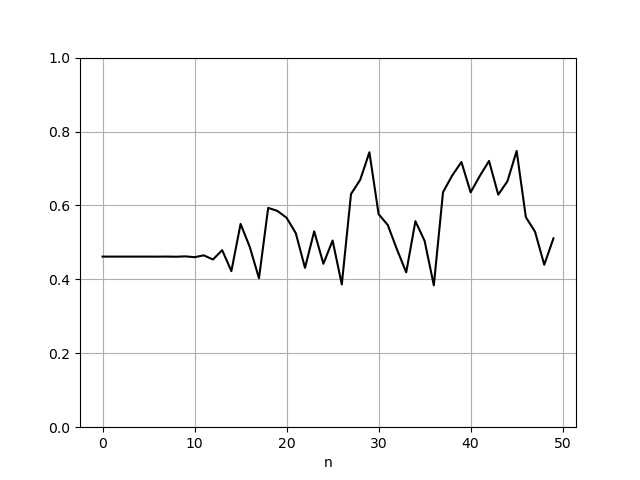}
        \caption*{$(T^2)^n(x_0)$}
        \label{fig:image2}
    \end{minipage}
    \caption{Trajectory of $T^n(x_0)$ and $(T^2)^n(x_0)$ for $x_0= 0.461538$}
    \label{fig:two_images}
\end{figure}

\FloatBarrier

To stabilize the unstable cycle of length \(2\) for the map \(T\), we stabilize the unstable fixed point $6/13$ of the map \(T^2\).

The stabilization algorithm is an A-stabilization method with recently proposed coefficients \cite{s}, based on the properties of a certain family of complex polynomials.

Let us describe the stabilization algorithm. We choose \(\sigma = 1.2\) and coefficients 
\begin{align*}
a_1 &= 6 c (\sigma^7 - \sigma^5)= 0.1854829091429525 \\
a_2 &= 5 c (3\sigma^7 - 5\sigma^5 + 2\sigma^3) = 0.2490279798678529 \\
a_3 &= 4 c (5\sigma^7 - 10 \sigma^5 + 6  \sigma^3 - \sigma)= 0.2485509147723208 \\
a_4 &= 3 c (5\sigma^7 - 10 \sigma^5 + 6  \sigma^3 - \sigma)= 0.1864131860792406 \\
a_5 &= 2 c (3\sigma^7 - 5  \sigma^5 + 2  \sigma^3)= 0.09961119194714116  \\
a_6 &= c(\sigma^7 - \sigma^5)= 0.030913818190492087 \\
\end{align*}
where constant $c$ was chosen to satisfy the condition that sum of all coefficients equal 1.
We claim that proposed coefficients have not approximately value but exactly equal to the writing numbers and one can easily check that indeed
\begin{equation}
\label{a_sum_1}
a_1+a_2+a_3+a_4+a_5+a_6=1
\end{equation}

Let us denote $f(x)=T^2(x)$ for simplicity. Now we start from the initial point $x_0\in[0,1]$ and then calculate first five points of the trajectory \(x_1=f(x_0)\), \(x_2=f^2(x_0)\), \(x_3=f^3(x_0)\), \(x_4=f^4(x_0)\), \(x_5=f^5(x_0)\). Then start the stabilization process.

 The next point will be a new point $x_6^*$ calculating by the rule
 $$
 x_n^*=a_1f(x_{n-1}^*)+a_2f(x_{n-2}^*)+ \cdots +a_6f(x_{n-6}^*)
 $$
for $n\geq 6$, where $x_n^*=x_n$ for $n=0,1,2,\dots 5$. We shall write $x^*_n=A_nf(x_0)$ as a result of stabilization.

Our aim was to study the global convergence of this stabilization algorithm. We hope that except some negligible set of points the process converges to one of the two cycle points. 

In our computer experiment, we chose 100,000 different initial points uniformly distributed over the interval $[0,1]$. Here are the results after 50 iterations for the initial points 0.1, 0.2, 0.3, and 0.4. These plots describe the general picture very well.

\begin{figure}[htbp]
    \centering
    \begin{minipage}[t]{0.48\textwidth}
        \centering
        \includegraphics[width=\textwidth]{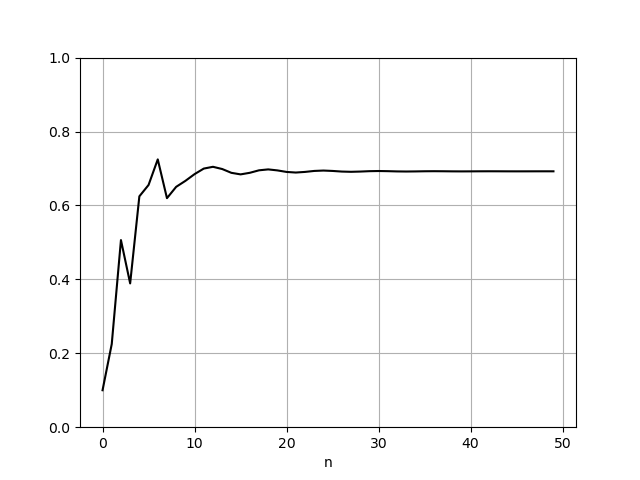}
        \caption*{$x_0=0.1$, $x_{50}^* \approx 0.6923$}
        \label{fig:image1}
    \end{minipage}
    \hfill
    \begin{minipage}[t]{0.48\textwidth}
        \centering
        \includegraphics[width=\textwidth]{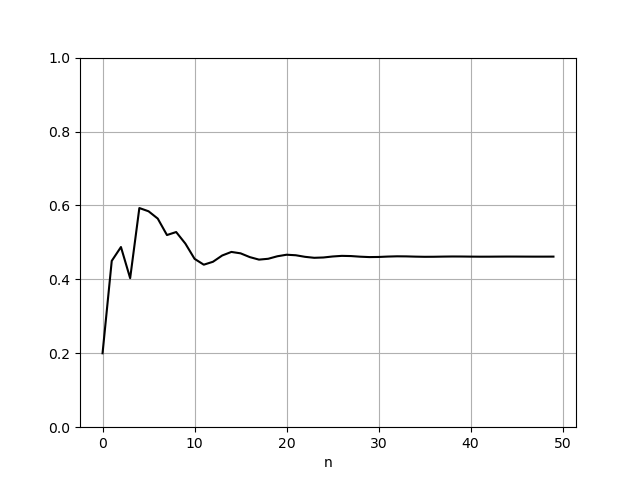}
        \caption*{$x_0=0.2$ , $x_{50}^* \approx 0.4615$  }
        \label{fig:image2}
    \end{minipage}
    \vspace{0.4cm}
    \begin{minipage}[t]{0.48\textwidth}
        \centering
        \includegraphics[width=\textwidth]{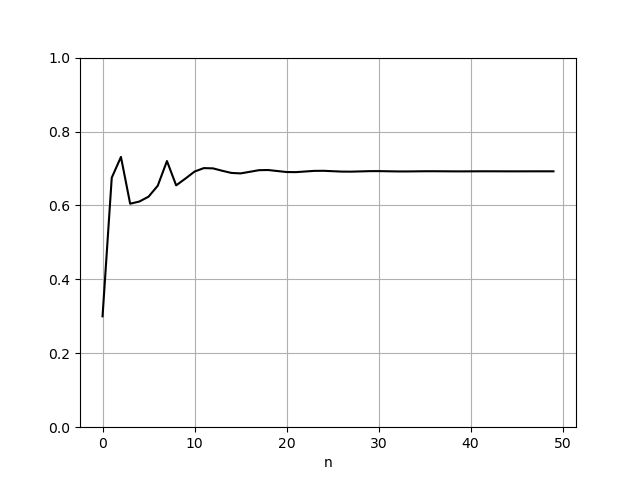}
        \caption*{$x_0=0.3$ , $x_{50}^* \approx 0.6923$ }
        \label{fig:image3}
    \end{minipage}
    \hfill
    \begin{minipage}[t]{0.48\textwidth}
        \centering
        \includegraphics[width=\textwidth]{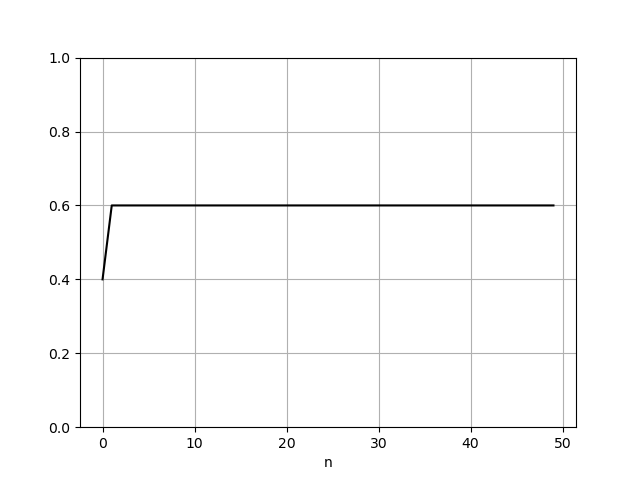}
        \caption*{$x_0=0.4$ , $x_{50}^* = 0.6$ }
        \label{fig:image4}
    \end{minipage}
    \caption{Behavior of the system $(T^2)^n$ after stabilization }
    \label{fig:four_images}
\end{figure}
\FloatBarrier 

Next are two plots: the first shows the point \( x^*_{50}=A_{50}f(x_0) \) against \( 10^5 \) initial points \( x_0 \), equally distributed over the segment \([0,1]\), after applying our stabilization algorithm. The second plot is a magnified view around one of the cycle points, \( 6/13=0.\overline{461538} \).

\begin{figure}[htbp]
    \centering
    \begin{minipage}[t]{0.48\textwidth}
        \centering
        \includegraphics[width=\textwidth]{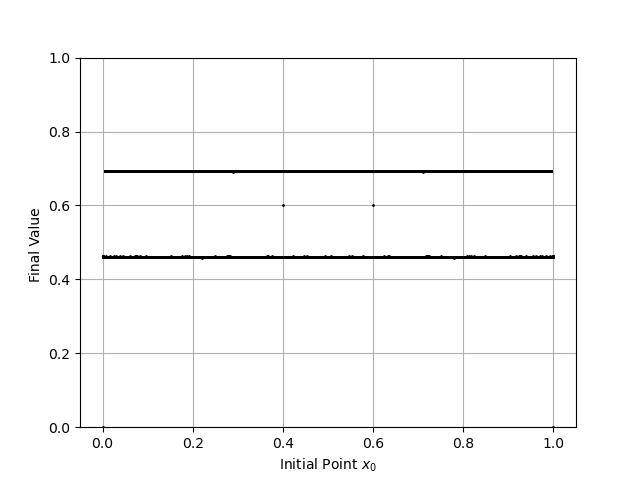}
        \caption*{$x^*_{50}$ against different $x_0$}
    \end{minipage}
    \hfill
    \begin{minipage}[t]{0.48\textwidth}
        \centering
        \includegraphics[width=\textwidth]{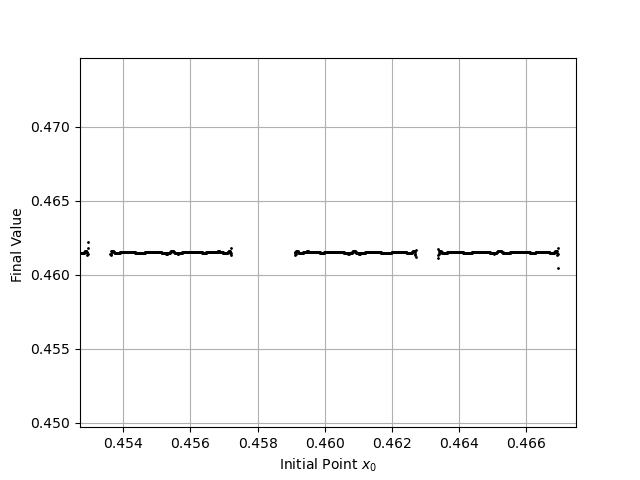}
        \caption*{$x^*_{50}$ for $x_0$ near $0.\overline{461538}$}
        \label{fig:image2}
    \end{minipage}
    \caption{Distribution of $x_{50}^*$ for different initial points  $x_0$}
    \label{fig:two_images}
\end{figure}

\FloatBarrier 
As one can clearly see, except for only $2$ points, namely $0.4$ and $0.6$, for all other initial conditions $x_0$, after just $50$ steps, the stabilizing sequence $A_n f(x_0)$ is approximately equal to one of the two points of the cycle. Actually, it can be verified that the difference between the calculated $x^*_{50}$ and the actual value of one of the points of the cycle is less than $10^{-3}$.

The behavior of the stabilization process starting at these $2$ exceptional points, $0.4$ and $0.6$, is very clear. Indeed, since $0.4$, after the first application of $f$, arrives at $0.6$, and $0.6$ is a fixed point of $f$, then all $x_1$, $x_2$, $x_3$, $x_4$, $x_5$, and $x_6$ are also equal to $0.6$. Hence, the average ( because of (\ref{a_sum_1})), also equal $0.6$. Thus, all points $x_n^*$ exactly equal $0.6$. The same is true for the initial point $x_0=0.6$.

Wait, but what about points like $x_0 = 4/15$? This point eventually falls to the fixed point of $T$, which equals $0.6$. Indeed, one can easily calculate that $x_1 = f(x_0) = T(T(x_0)) = T(6/10) = 6/10$, $x_2 =f(x_1)= 6/10$, and so on. Hence, the averaging process is just equal to $0.6$ at every step. Why did we not see this point even after choosing $10^{10}$ equidistributed initial points? Moreover, there are many such points, like $11/15$ or $23/45$, which eventually become fixed for $T$. Thus, our averaging process must converge to the fixed point of $T$ equal $0.6$.

The explanation is clear: points like $4/15 = 0.266666\ldots$ are not included in any of our equidistributed nets of $10^m$ initial points!

Let us choose $5\cdot3^m$ equidistributed initial points for $m = 5$ or $m = 10$ and examine the distribution of $x^*_{30}$ (30 iterations for more clear picture).

\begin{figure}[htbp]
    \centering
    \begin{minipage}[t]{0.48\textwidth}
        \centering
        \includegraphics[width=\textwidth]{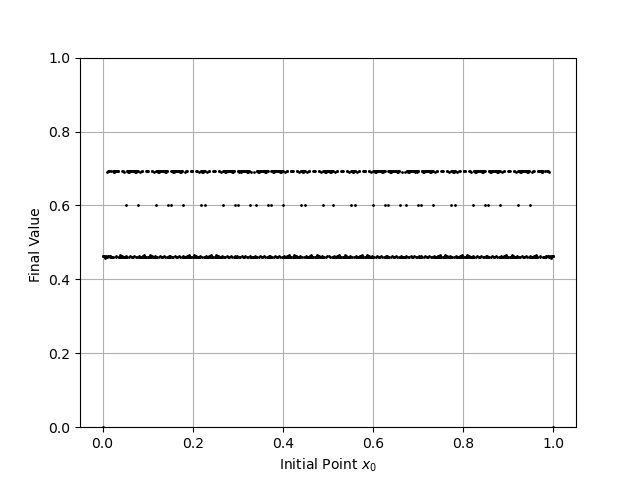}
        \caption*{ $x^*_{30}$ for $5\cdot3^{5}$ points $x_0$}
        \label{fig:image1}
    \end{minipage}
    \hfill
    \begin{minipage}[t]{0.48\textwidth}
        \centering
        \includegraphics[width=\textwidth]{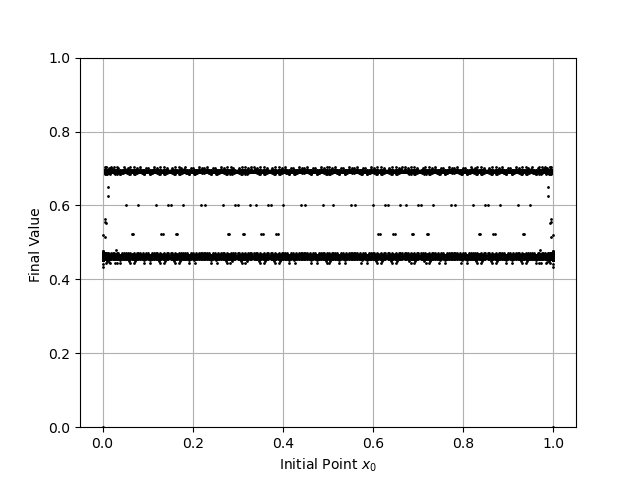}
        \caption*{$x^*_{30}$ for $5\cdot3^{10}$ points $x_0$}
        \label{fig:image2}
    \end{minipage}
    \caption{Distribution of $x_{30}^*$ for equidistributed $5\cdot3^m$ initial points  $x_0$}
    \label{fig:two_images}
\end{figure}

\FloatBarrier 

The conclusions from these experiments allow us to formulate a working hypothesis that, for all initial points except for a certain countable set, the stabilization sequence converges to one of the cycles.

To increase the accuracy of calculating the cycle points, the author increased the number of iterations of the stabilization algorithm to $100$. 

We will present only two graphs: one corresponding to the initial point \(0.3\) and the other to the initial point \(0.4\). While the first graph raises no questions (the accuracy improved, and \(x_{100} =0.692307 \pm 10^{-6}\)), something strange happened with the second graph.

\begin{figure}[htbp]
    \centering
    \begin{minipage}[t]{0.48\textwidth}
        \centering
        \includegraphics[width=\textwidth]{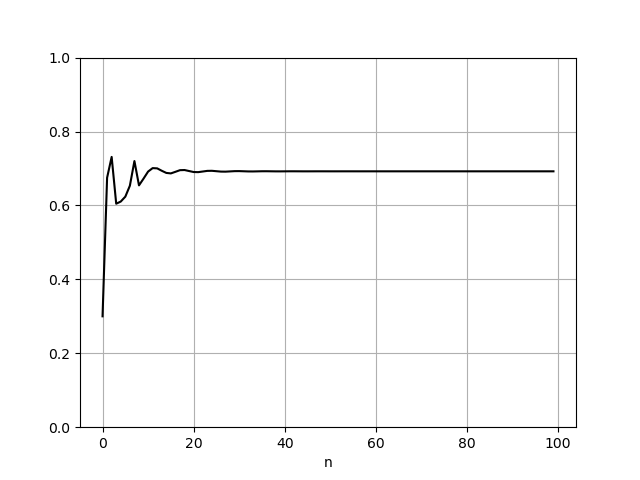}
        \caption*{Initial point $x_0=0.3$}
        \label{fig:image1}
    \end{minipage}
    \hfill
    \begin{minipage}[t]{0.48\textwidth}
        \centering
        \includegraphics[width=\textwidth]{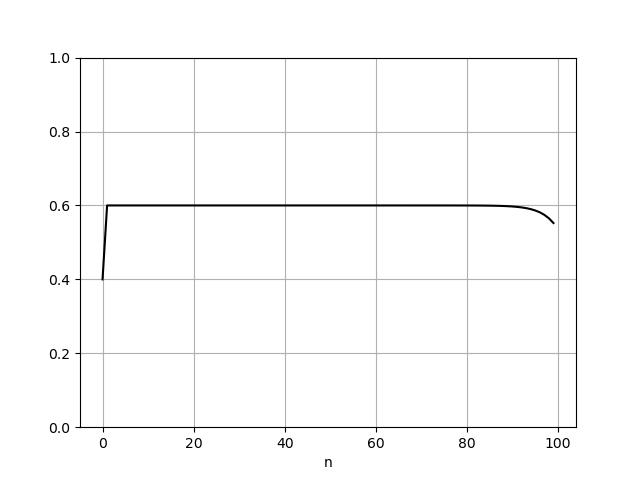}
        \caption*{Initial point $x_0=0.4$}
        \label{fig:image2}
    \end{minipage}
    \caption{Two plots of stabilization starting at points $0.3$ and $0.4$}
    \label{fig:two_images}
\end{figure}

\FloatBarrier

Instead of remaining a constant function indefinitely, we observe something strange starting from the $90$th step.

\section{The Rabbit Hole}

Then, we increased the number of iterations to 150. Oops, what is happening?

\begin{figure}[htbp]
    \centering
    \includegraphics[width=\textwidth]{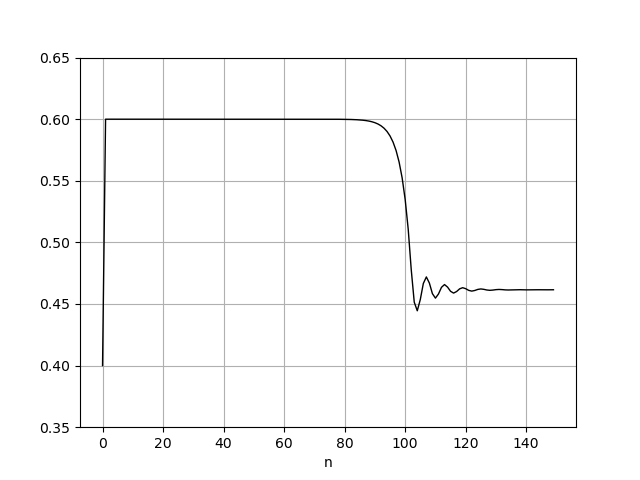} 
    \caption{Plot of $A_nf(0.4)$ for $n\leq 150$}
    \label{fig:plot}
\end{figure}

\FloatBarrier

\textit{The stabilized system, after remaining fixed for approximately 100 units of time, suddenly transitions to one of its stable points.}

There are two pieces of news after the discovery of this phenomenon—one bad and one good.

The bad news is that when observing a fixed state or cycle after stabilization, we can’t guarantee that this state or cycle will remain stable for a significant amount of time.

But there is also good news. According to proven theorems, the stabilization algorithm allows us to detect cycles for any initial values, except for a countable (but infinite) set of states. Thus, this holds true for almost all initial conditions.

The good news is that, \textit{in practice}, the algorithm \textit{always} finds cycles.

Now, we shall provide a sketch of the proof that the algorithm converges except for some negligible set of initial points and also explain the observed phenomenon.

Let us define
\begin{equation}
\label{pro_a_stab_system}
F(u_1, u_2, \dots, u_6) = \big( u_2, u_3, \dots, u_6, a_1 f(u_6) + a_2 f(u_5) + \dots + a_6 f(u_1) \big).
\end{equation}

Then, \(F\) maps the 6-dimensional cube \(Q = [0,1]^6\) to itself. It is clear that if \(s\) is a fixed point of \(f\), then \((s, s, \dots, s)\) is a fixed point of \(F\). Conversely, if \((s_1, s_2, s_3, \dots, s_6)\) is a fixed point of \(F\), then \(s_1 = s_2 = \dots = s_6 = s\), and \(s\) is a fixed point of \(f\).

Thus, the discrete dynamical system \( F: Q \to Q \) has only three fixed points, namely \( p = (6/13, \dots, 6/13) \), \( q = (9/13, \dots, 9/13) \), and \( r = (3/5, \dots, 3/5) \). Let us divide \( Q \) into \( 2^6 = 64 \) equal boxes \( Q_i \), within which \( F \) acts as an affine map (linear plus constant). Two of these boxes, \( Q_p \) and \( Q_q \), contain \( p \) and both \( q \) and \( r \), respectively. Since \( F \) is essentially a linear operator, its dynamics are governed by its eigenvalues. However, we chose the coefficients \(a_k\) in such a manner that \( p \) and \( q \) become stable points, while \( r \) does not. 

Sooner or later (except for a negligible set of initial points), the trajectory of an initial point will reach one of these boxes \( Q_p \) or \( Q_q \) and be attracted to one of the stable points.

To explain the observed phenomenon, we note that due to precision limitations (approximation), the trajectory does not fall exactly onto the point \( r \) but instead into an extremely small neighborhood around it. Meanwhile, because \( r \) is unstable, the trajectory eventually leaves this neighborhood and is attracted to one of the stable points.

As a one more example, consider the graph of the orbit \( T^n_h(0.5) \) when \( h \) equals the square root of two, calculated to 50 decimal places of precision:
\begin{equation}
\label{appr_sqrt2}
h = 1.414213562373095048801688724209698078569671875376948073176
\end{equation}

The left graph shows the orbit after the first 300 iterations, clearly demonstrating that the point \( 0.5 \) eventually becomes fixed. The right graph extends this by showing the orbit over the next 300 iterations of the tent map.

\begin{figure}[htbp]
    \centering
    \begin{minipage}[t]{0.48\textwidth}
        \centering
        \includegraphics[width=\textwidth]{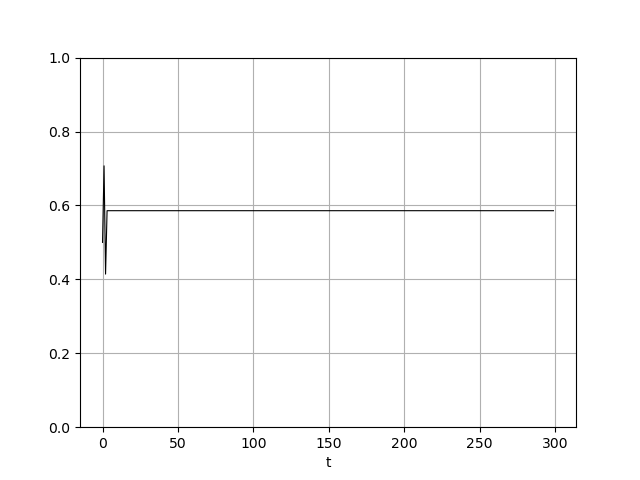}
        \caption*{First 300 iterations}
        \label{fig:image1}
    \end{minipage}
    \hfill
    \begin{minipage}[t]{0.48\textwidth}
        \centering
        \includegraphics[width=\textwidth]{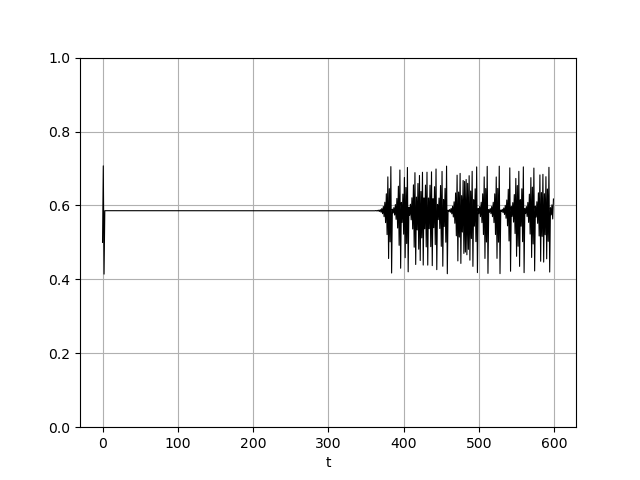}
        \caption*{Next 300 iterations}
        \label{fig:image2}
    \end{minipage}
    \caption{Plot of $T^n_h(0.5)$ }
    \label{fig:two_images}
\end{figure}

\FloatBarrier 

Here not only dynamics of the "critical point" $x_0=0.5$ change significantly after we pass from $\sqrt{2}$ to the $h$ defined as (\ref{appr_sqrt2}), but also all dynamical picture dramatically change. This is a bifurcation point of the family $T_h$.

Finally we would like to explain the title of the article.

\section{Fibonacci Rabbits}
The numbers in the sequence  
\[
1, 1, 2, 3, 5, 8, 13, 21, 34, \dots
\]  
are called Fibonacci numbers, in which each number is the sum of the preceding two. They can be formally defined as:  
\begin{equation}
\label{fib}
x_0 = 1, \quad x_1 = 1, \quad x_n = x_{n-1} + x_{n-2} \quad \text{for} \quad n \geq 2.
\end{equation}

This famous sequence was published in 1202 by Leonardo Pisano (Leonardo of Pisa), who is sometimes called Leonardo Fibonacci (Filius Bonaccii, "son of Bonaccio"). His \textit{Liber Abaci} (\textit{The Book of Calculation}) contains the following exercise: "How many pairs of rabbits can be produced from a single pair in a year's time?"  

To solve this problem, we are told to assume that each pair produces a new pair of offspring every month and that each new pair becomes fertile at the age of one month. Furthermore, the rabbits never die. After one month, there will be 2 pairs of rabbits; after two months, there will be 3. The following month, the original pair and the pair born during the first month will each usher in a new pair, resulting in 5 pairs in total, and so on.

These are the  White Rabbits. 

The Black Rabbits appear when we choose \( x_0 = 1 \) and \( x_1 = -0.618033988749 \) in (\ref{fib}). Here $x_1$ is an approximation of $-1/\phi$, where $\phi$ is the golden ratio.

\begin{figure}[htbp]
    \centering
    \begin{minipage}[t]{0.48\textwidth}
        \centering
        \includegraphics[width=\textwidth]{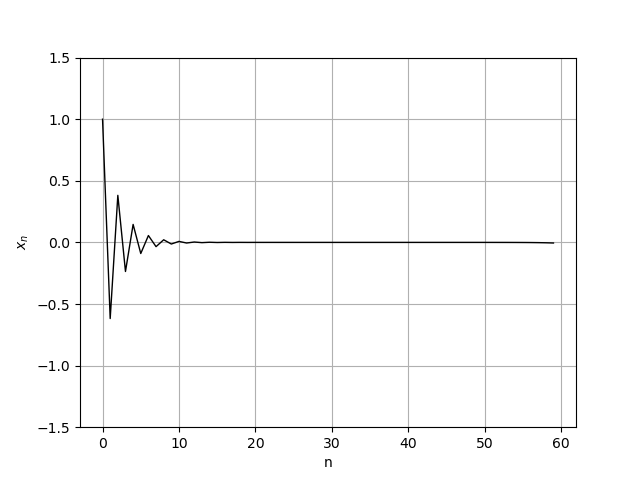}
        \caption*{$x_n$ before $60$}
        \label{fig:image1}
    \end{minipage}
    \hfill
    \begin{minipage}[t]{0.48\textwidth}
        \centering
        \includegraphics[width=\textwidth]{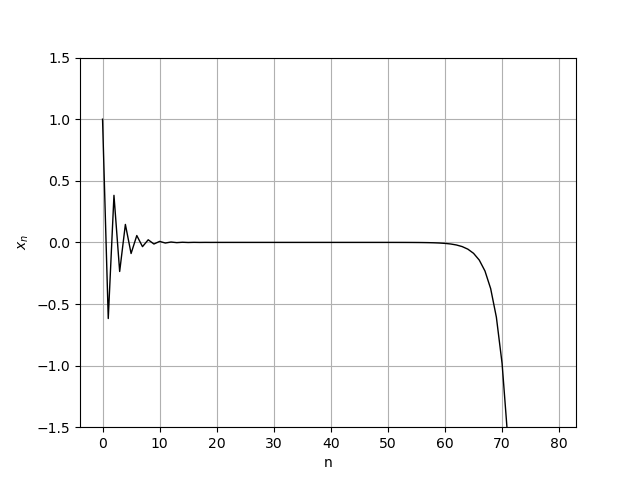}
        \caption*{$x_n$ after $60$}
        \label{fig:image2}
    \end{minipage}
    \caption{Black Rabbits of Fibonacci}
    \label{fig_brf}
\end{figure}

Such defined sequence demonstrate a very stable and predictable behavior (left plot in Fig. \ref{fig_brf} ). It clearly tends to zero and only after $40$ steps became less that $0.00001$ in absolute value

\FloatBarrier 
But after $n=60$ the sequence suddenly fall down.

In this model case, it is easy to provide an explanation for this phenomenon. It is sufficient to transition to a higher-dimensional phase space and consider the behavior of the points $(x_n, x_{n+1})$. These points lie on a cross formed by two intersecting lines.

\begin{figure}[htbp]
    \centering
    \includegraphics[width=\textwidth]{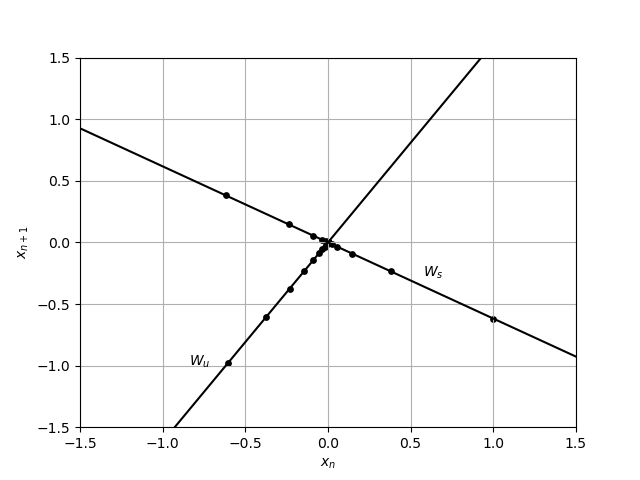} 
    \caption{Points $(x_n, x_{n+1})$ and stable ($W_s$) and unstable ($W_u$) manifolds }
    \label{fig:plot}
\end{figure}

\FloatBarrier

These lines serve as the stable manifold ($W_s$) and the unstable manifold ($W_u$) for the linear operator
\[
A=
\begin{pmatrix}
0 & 1 \\
1 & 1
\end{pmatrix}
\]
which governs the dynamics of the sequence $x_n$ in the sense that
\[
A
\begin{pmatrix} 
x_{n-2} \\ 
x_{n-1} 
\end{pmatrix}
= \begin{pmatrix} 
x_{n-1} \\ 
x_n 
\end{pmatrix}.
\]

This linear operator has two eigenvectors corresponding to the lines \(W_u\) and \(W_s\), and two eigenvalues, namely the golden ratio \(\phi\) and \(-1/\phi\). One of these eigenvalues has an absolute value greater than one, while the other has an absolute value less than one. Such operators are called hyperbolic.

When a point \((x_n, x_{n+1})\) is close to the stable manifold, it is attracted to zero, since \( \lvert -1/\phi \rvert < 1 \). However, as soon as it approaches the unstable manifold, the dynamics shift, and the point is pulled away to infinity.

And finally, a bit of philosophy.

In his book \textit{The Black Swan: The Impact of the Highly Improbable}, Nassim Nicholas Taleb explores the profound and unpredictable events that shape our world, which he calls ``Black Swans.'' These events are characterized by three main features: they are rare and unpredictable, they have an extreme impact, and in hindsight, they often seem explainable or inevitable.

Taleb argues that our world is increasingly shaped by these Black Swans, from financial crises to technological breakthroughs, and yet humans are inherently ill-equipped to anticipate them due to cognitive biases and reliance on standard predictive models.

By analogy with this term, we have named Black Rabbits events that exhibit all the external characteristics of a Black Swan but are entirely determined by the state of the system itself.

These phenomena are not artificial examples or mere curiosities; they exist wherever hyperbolic meadows are found.

When we witness a social catastrophe or a financial crisis, could it be not a stray Black Swan but a local Black Rabbit?

How can one determine in advance that a system is pregnant with a crisis? What indicators can be used? Is it possible to control the system to prevent it from becoming unstable?

These are intriguing questions, perhaps too complex, especially when dealing not with simple mathematical models but with real-world systems. Yet, we must still strive to find answers to them.


\begin{thebibliography}{9}

\bibitem{cg}
Lennart Carleson and Theodore W. Gamelin, 
\textit{Complex Dynamics}, 
Springer, Universitext Series, New York, 1993. DOI: \url{10.1007/978-1-4419-3117-6}.

\bibitem{ce}
Pierre Collet and Jean-Pierre Eckmann, 
\textit{Iterated Maps of the Interval as Dynamical Systems}, 
Birkhäuser, Progress in Physics Series, Volume 1, Boston, MA, 1980. DOI: \url{10.1007/978-1-4612-5929-5}.


\bibitem{mt}
John Milnor and William Thurston, 
\textit{On Iterated Maps of the Interval}, 
Lecture Notes in Mathematics, Springer-Verlag, Volume 1342, 1988, pp. 465–563. DOI: \url{https://doi.org/10.1007/BFb0082847}.



\bibitem{s}
Alexei Solyanik, 
\textit{M-Stable Polynomials}, 
Plenary talk at the International Conference on Approximation and Potential Theory, 
Georgia Southern University, September 22--24, 2023. 
\textit{Conference Website:} \url{https://sites.google.com/georgiasouthern.edu/icapt2023/home}.

\end{thebibliography}
\end{document}